\documentclass[12pt]{article}
\usepackage{amssymb, amsmath, amsfonts, tikz, geometry, ytableau}
\usepackage[indent,margin=1cm]{caption}
\usepackage[colorlinks,linkcolor=blue,anchorcolor=blue,citecolor=blue]{hyperref}

\numberwithin{equation}{section}
\numberwithin{figure}{section}

\hypersetup{colorlinks = true}
\newtheorem{thm}{Theorem}[section]
\newtheorem{conj}[thm]{Conjecture}

\newtheorem{cor}[thm]{Corollary}
\newtheorem{lem}[thm]{Lemma}
\newtheorem{exam}[thm]{Example}
\newtheorem{prop}[thm]{Proposition}

\newtheorem{pf}[thm]{Proof}

\def\pf{\noindent{\it Proof.} }
\def\qed{\nopagebreak\hfill{\rule{4pt}{7pt}}
	\medbreak}

\allowdisplaybreaks

\def\pf{\noindent{\it Proof.} }
\def\qed{\nopagebreak\hfill{\rule{4pt}{7pt}}
\medbreak}

\allowdisplaybreaks

\begin{document}
\begin{center}
	{\large \bf Immanant Positivity for Catalan-Stieltjes Matrices}
\end{center}

\begin{center}
Ethan Y.H. Li$^{1}$, Grace M.X. Li$^{2}$,  Arthur L.B. Yang$^{3}$ and Candice X.T. Zhang$^{4}$\\[6pt]

Center for Combinatorics, LPMC\\
Nankai University, Tianjin 300071, P. R. China\\[8pt]

Email: $^{1}${\tt yinhao\_li@mail.nankai.edu.cn}, $^{2}${\tt limengxing@mail.nankai.edu.cn}, $^3${\tt yang@nankai.edu.cn}, $^{4}${\tt zhang\_xutong@mail.nankai.edu.cn}
\end{center}

\noindent\textbf{Abstract.}
In this paper we give some sufficient conditions for
the nonnegativity of immanants of square submatrices of Catalan-Stieltjes matrices and their corresponding Hankel matrices.
To obtain these sufficient conditions, we construct new planar networks with a recursive nature for Catalan-Stieltjes matrices.
As applications, we provide a unified way to produce inequalities for many combinatorial polynomials, such as the Eulerian polynomials, Schr\"oder polynomials and Narayana polynomials.

\noindent \emph{AMS Mathematics Subject Classification 2020:} 05A05, 05A20

\noindent \emph{Keywords:}  immanant, character, Catalan-Stieltjes matrices, Hankel matrices, planar network

\section{Introduction}

In recent years the ($q$-)total positivity of Catalan-Stieltjes matrices and the associated Hankel matrices has been extensively studied, see \cite{Chen-recursive-matrices-2015,Pan-Zeng-2016, Wang-Zhu-2016} and the references therein. The main objective of this paper is to study the
($q$-)nonnegativity of immanants of these matrices, which is motivated by Stembridge's result on the nonnegativity of immanants of totally positive matrices \cite{Stembridge-1991}. Before stating our main result, let us first give an overview of some related concepts and results.

A \emph{Catalan-Stieltjes matrix} is an infinite lower triangular matrix $(c_{n,k})_{n\ge k \ge 0}$, which can be recursively defined by
\begin{equation}\label{recurrence-CS}
\begin{split}
c_{n,0} &=s_0 c_{n-1,0}+t_1c_{n-1,1};\\
c_{n,k}&=r_{k-1}c_{n-1,k-1}+s_kc_{n-1,k}+t_{k+1}c_{n-1,k+1} \quad (k\ge 1, \, n\ge 1),
\end{split}
\end{equation}
where $c_{0,0}=1$ and $r_k,s_k,t_{k+1}$ are certain given parameters for $k\geq 0$.
Letting $\gamma=(r_k)_{k\ge 0}$, $\sigma=(s_k)_{k\ge 0}$ and $\tau=(t_k)_{k\ge 1}$, we usually use $C^{\gamma,\sigma,\tau}$ to denote the matrix $(c_{n,k})_{n\ge k \ge 0}$.
The term ``Catalan-Stieltjes matrix'' was coined by Pan and Zeng \cite{Pan-Zeng-2016}
when $\gamma$, $\sigma$ and $\tau$ consist of nonnegative integers, while Aigner called $C^{\gamma,\sigma,\tau}$ a \emph{Catalan matrix} in \cite{Aigner-a-course}, a \emph{recursive matrix} in \cite{Aigner-2001}, and called the numbers $c_{n,0}\, (n\geq 0)$   \emph{Catalan-like numbers} in \cite{Aigner-SM-1999}.
If $\gamma$, $\sigma$ and $\tau$ are three sequences of polynomials in one variable, say $q$, Wang and Zhu \cite{Wang-Zhu-2016} called $C^{\gamma,\sigma,\tau}$ a \emph{$q$-recursive matrix} and the entries $c_{n,0}\,(n\geq 0)$  \emph{$q$-Catalan-like numbers}. By abuse of notation, in this paper we will follow Pan and Zeng to call
$C^{\gamma,\sigma,\tau}$ of the form \eqref{recurrence-CS} a \emph{Catalan-Stieltjes matrix} and $c_{n,0}\,(n\geq 0)$ \emph{Catalan-like numbers} no matter whether $\gamma$, $\sigma$ and $\tau$ depend on some indeterminates or not. 
We use $H^{\gamma,\sigma,\tau}$ to denote the associated Hankel matrix $(c_{i+j,0})_{i,j \ge 0}$.

The study of totally positive matrices arises in various branches of mathematics, probability, statistics, mechanics, economics and computer science, for more information see \cite{Karlin-1968}.
Recall that a matrix is said to be \emph{totally positive}, or TP for short, if all of its minors are nonnegative numbers, and it is said to be \emph{$q$-totally positive}, or $q$-TP for short, if each of its minors is a polynomial in $q$ with nonnegative coefficients.
For convenience, if $f(q)$ is a polynomial in $q$ with nonnegative coefficients, then we say that it is $q$-\emph{nonnegative}. Moreover, we write $f(q)\geq_q g(q)$ if $f(q)-g(q)$ is $q$-nonnegative.
A sequence $\alpha=(a_k)_{k\ge 0}$ consisting of nonnegative numbers (respectively, $q$-nonnegative polynomials) is said to be \emph{log-convex} (respectively, $q$-\emph{log-convex}) if $a_ka_{k+2}\ge a_{k+1}^2$ (respectively, $a_ka_{k+2}\ge_q a_{k+1}^2$) for all $k\ge 0$.
A sequence $\alpha=(a_k)_{k\ge 0}$  is called \emph{Hankel totally positive} (respectively, $q$-\emph{Hankel totally positive}) if all the minors of its Hankel matrix $(a_{i+j})_{i,j\geq 0}$ are nonnegative (respectively, $q$-nonnegative).
We usually abbreviate them as H-TP and $q$-H-TP, respectively.
Note that $q$-totally positive matrices are also said to be  coefficientwise totally positive in $q$, and $q$-{Hankel totally positive} sequences are also said to be
coefficientwise Hankel totally positive in $q$, see \cite{Chen-coef-2021,Sokal-1980}.
A sequence $\alpha=(a_k)_{k\ge 0}$ is said to be a Stieltjes moment (respectively, $q$-Stieltjes moment) sequence if it is H-TP (respectively, $q$-H-TP).
It is clear that a Stieltjes moment sequence (respectively, $q$-Stieltjes moment sequence) must be log-convex (respectively, $q$-log-convex).
We would like to point out that
Stieltjes moment sequences and $q$-Stieltjes moment sequences arise in many fields of mathematics, see for instance \cite{Aigner-2001,Bennett-SM-2011,
Catalan-like-2016,Shohat-moments-1943, Widder-laplace-1946}.

Since numbers can be considered as  polynomials of degree zero, the notion of total positivity can be treated as a special case of that of $q$-total positivity. For this reason, we only state our results on $q$-total positivity throughout this paper, unless specifically declared. The following sufficient conditions for the $q$-total positivity of Catalan-Stieltjes matrices and the corresponding Hankel matrices have been obtained.

\begin{thm}[{\cite[Corollary 2.4]{Pan-Zeng-2016},\cite[Lemma 3.3]{Wang-Zhu-2016}}]
\label{thm-tp-csmatrix}
If the sequences $\gamma=(r_k)_{k\ge 0}$, $\sigma=(s_k)_{k \ge 0}$ and $\tau=(t_k)_{k\ge 1}$ of $q$-nonnegative polynomials satisfy one of the following conditions:
  \begin{itemize}
	\item[(1)] $s_0\ge_q r_0$, and $s_k\ge_q r_k+t_k$ for $k\ge 1$;
	
	\item[(2)] $s_0\ge_q t_1$, and $s_k\ge_q r_{k-1}+t_{k+1}$ for $k\ge 1$;
	
	\item[(3)] $s_0\ge_q 1$, and $s_k\ge_q r_{k-1}\cdot t_{k}+1$ for $k\ge 1$;
	
	\item[(4)] $s_0\ge_q r_0\cdot t_1$, and $s_k\ge_q r_k\cdot t_{k+1}+1$ for $k\ge 1$;

	\item[(5)] there exist two $q$-nonnegative polynomials $b_k$ and $c_k$ such that $r_k=1$, $s_k=b_k+c_k$, and $t_{k+1}=b_{k+1}c_k$ for each $k\ge 0$,
  \end{itemize}
  then both $C^{\gamma,\sigma,\tau}$ and $H^{\gamma,\sigma,\tau}$ are $q$-TP.
\end{thm}

Liang, Mu and Wang \cite{Catalan-like-2016} first proved that a sequence of ordinary Catalan-like numbers is H-TP under the first and third conditions with $r_k = 1$ ($k \ge 0$).
Chen, Liang and Wang \cite{Chen-recursive-matrices-2015} established the $q$-total positivity of
$C^{\gamma,\sigma,\tau}$ in the first three cases of Theorem \ref{thm-tp-csmatrix}.
The $q$-total positivity of
$C^{\gamma,\sigma,\tau}$ in the fourth case and the $q$-total positivity of
$H^{\gamma,\sigma,\tau}$ in the first four cases were given by Pan and Zeng \cite{Pan-Zeng-2016}. Wang and Zhu \cite{Wang-Zhu-2016} first proved the $q$-total positivity of $H^{\gamma,\sigma,\tau}$ explicitly in the fifth case, as well as the $q$-total positivity of $C^{\gamma,\sigma,\tau}$ implicitly.

Now we turn to review
some results on immanants.
Recall that, given a matrix $M=(m_{i,j})_{1\le i,j \le n}$ and a partition $\lambda$ of $n$,
the immanant of $M$ with respect to $\lambda$ is defined by
\begin{equation*}
{\rm Imm}_\lambda M=\sum_{\pi\in\mathfrak{S}_n}\chi^{\lambda}(\pi)\prod_{i=1}^n m_{i,\pi(i)},
\end{equation*}
where $\pi$ ranges over all permutations in the symmetric group $\mathfrak{S}_n$ and $\chi^{\lambda}$ denotes the irreducible character of  $\mathfrak{S}_n$ associated with $\lambda$. The concept of an immanant was introduced by Littlewood \cite{Littlewood} to define the Schur symmetric functions in terms of the power sum symmetric functions. It is clear that the immanant ${\rm Imm}_{(1^n)} M$ specializes to $\det M$, the determinant of $M$.
Goulden and Jackson \cite{GJ92} initiated
the study of the positivity of immanants of combinatorial matrices. Stembridge \cite{Stembridge-1991} proved that any immanant of any totally positive matrix is nonnegative. Based on a planar network interpretation for
totally positive matrices, Brenti \cite{Brenti-1995} provided a different proof of this fact. Both proofs depend on a result due to Cryer \cite{Cryer1976},
which states that each TP matrix can be written as a product of TP bidiagonal matrices. However, Cryer's result fails for
$q$-TP matrices, more precisely, not every $q$-TP matrix can be written as a product of $q$-TP bidiagonal matrices,
and it is still unknown whether
the immanants of a $q$-totally positive matrix are $q$-nonnegative.

The main result of this paper is as follows, which obviously generalizes Theorem \ref{thm-tp-csmatrix}.

\begin{thm}\label{thm-main}
If the three sequences $\gamma=(r_k)_{k\ge 0}$, $\sigma=(s_k)_{k \ge 0}$ and $\tau=(t_k)_{k\ge 1}$ consist of $q$-nonnegative polynomials and satisfy one of the five conditions in Theorem \ref{thm-tp-csmatrix}, then every immanant of each square submatrix of $C^{\gamma,\sigma,\tau}$ is $q$-nonnegative. The same is true for every immanant of each square submatrix of $H^{\gamma,\sigma,\tau}$.
\end{thm}

\noindent \textbf{Remark.} The nonnegativity of the immanants of certain Hankel matrices $H^{\gamma,\sigma,\tau}$ was also studied by Goulden and Jackson \cite{GJ92}.

We would like to point out that
the key step of Brenti's proof of the immanant positivity of a TP matrix is to
construct a planar network for this matrix.
In order to prove Theorem \ref{thm-main},
it is certainly desirable to provide suitable planar networks for $C^{\gamma,\sigma,\tau}$ and  $H^{\gamma,\sigma,\tau}$. Actually, Chen, Liang and Wang \cite{Chen-recursive-matrices-2015} already raised the problem of finding a combinatorial interpretation for the total positivity of Catalan-Stieltjes matrices, and subsequently, Pan and Zeng \cite{Pan-Zeng-2016} provided a planar network interpretation for the first four cases of
Theorem \ref{thm-tp-csmatrix}. Here we will
give a unified approach to Chen, Liang and Wang's problem for all cases of Theorem \ref{thm-tp-csmatrix}. Not only is our construction natural, but it can also be easily used to deduce Pan and Zeng's planar network interpretation.

Theorem \ref{thm-tp-csmatrix} can be used to prove the $q$-log-convexity or $q$-Stieltjes moment property of many combinatorial numbers or polynomials, such as Catalan numbers and Narayana polynomials, see
\cite{Wang-Zhu-2016, Zhu-2013} for instance.  Since Theorem \ref{thm-main} is stronger than Theorem \ref{thm-tp-csmatrix}, it is natural to expect that one can get more properties of related sequences and polynomials from the results on the immanants of their square submatrices.

The paper is organized as follows.
Section \ref{sect-ci} will be devoted to the recursive construction of planar networks for Catalan-Stieltjes matrices and their Hankel matrices.
In Section \ref{sect-im} we will
prove Theorem \ref{thm-main} based on our combinatorial interpretation of these matrices, and then apply Theorem \ref{thm-main} and related results to some combinatorial polynomials to obtain new inequalities.

\section{Planar networks for \texorpdfstring{$C^{\gamma,\sigma,\tau}$}{} and \texorpdfstring{$H^{\gamma,\sigma,\tau}$}{}} \label{sect-ci}

The aim of this section is to give a combinatorial interpretation for $C^{\gamma,\sigma,\tau}$ and $H^{\gamma,\sigma,\tau}$ with respect to the conditions of Theorem \ref{thm-tp-csmatrix}. Specifically, we construct a planar network for the Catalan-Stieltjes matrix $C^{\gamma,\sigma,\tau}$ with $q$-nonnegative weights for each case of Theorem \ref{thm-tp-csmatrix}.
With a small modification, our approach can be applied to construct a planar network for  $H^{\gamma,\sigma,\tau}$ with $q$-nonnegative weights. For the case of
$\gamma$ consisting of only $1$'s, we provide an alternative planar network for the Hankel matrix $H^{\gamma,\sigma,\tau}$.

In order to describe the planar networks given here, let us first introduce some definitions. By a planar network we mean
a quadruple $\mathcal{D}=\left(D, \mathrm{wt}_D,\mathbf{U},\mathbf{V}\right)$, where $D$ is an acyclic planar and locally finite digraph with vertex set $V(D)$ and arc set $A(D)$, $\mathrm{wt}_D$, called a weight function of $D$, is a map from $A(D)$ to a commutative ring with identity, and $\mathbf{U}$ and $\mathbf{V}$ are two sequences of vertices in $D$. Given two vertices $u,v$ of $D$, we denote by $GF_\mathcal{D}(u,v)$ the sum of the weights of all directed paths from $u$ to $v$, where the weight of a path is the product of the weights of all its arcs. In particular, we put $GF_{\mathcal{D}}(u,u) = 1$ by convention.
Given an $l\times m$ matrix $X=(x_{i,j})$, if there exists a planar network
$\mathcal{X}=(D^X, \mathrm{wt}_{D^X}, (u_1,u_2,\ldots,u_l),(v_1,v_2,\ldots,v_m))$
such that
\begin{equation*}
  x_{i,j} = GF_{\mathcal{X}}(u_i,v_j),
\end{equation*}
then $\mathcal{X}$ is called a planar network for $X$.
By the transfer-matrix method (see \cite[Theorem 4.7.1]{StaEC1} for instance), we have the following result.

\begin{lem}\label{lem-transfer-matrix}
Given an $l\times m$ matrix $X$ and an $m\times n$ matrix $X'$, suppose that
$\mathcal{X}=(D^X, \mathrm{wt}_{D^X}, \mathbf{U},\mathbf{V})$
is a planar network for $X$ and $\mathcal{X'}=(D^{X'}, \mathrm{wt}_{D^{X'}}, \mathbf{U'},\mathbf{V'})$
is a planar network for $X'$, where $\mathbf{U} = (u_1,u_2,\ldots,u_l)$, $\mathbf{U'} = (u'_1,u'_2,\ldots,u'_m)$ are two sequences of sources in $D^X$ and $D^{X'}$, respectively, and $\mathbf{V} = (v_1,v_2,\ldots,v_m)$, $\mathbf{V'} = (v'_1,v'_2,\ldots,v'_n)$ are two sequences of sinks in $D^X$ and $D^{X'}$, respectively.
Let $D^{XX'}$ be the union of $D^{X}$ and $D^{X'}$ with $v_i$ and $u_i'$ being identified for each $1\leq i\leq m$, and let $\mathrm{wt}_{D^{XX'}}$ be the weight function inherited from $\mathrm{wt}_{D^{X}}$ and $\mathrm{wt}_{D^{X'}}$.
Then
$(D^{XX'}, \mathrm{wt}_{D^{XX'}}, \mathbf{U}$, $\mathbf{V'})$
is a planar network for $XX'$.
\end{lem}

\pf It follows from the construction of $D^{XX'}$ and the choices of $\mathbf{U}$, $\mathbf{V}$, $\mathbf{U'}$ and $\mathbf{V'}$ that each directed path from $u_i$ to $v_j'$ must pass through one and only one $v_k(=u_k')$ for some $1 \le k \le m$. Hence
$$GF_{\mathcal{XX'}}(u_i,v_j') = \sum_{k=1}^m GF_{\mathcal{X}}(u_i,v_k)GF_{\mathcal{X'}}(u_k',v_j').$$ This completes the proof. \qed

Before describing our construction of planar networks, let us first review some properties of Catalan-Stieltjes matrices.
Let $C_n=(c_{i,j})_{0\leq i,j\leq n}$ be the $n$th leading principal submatrix of $C^{\gamma,\sigma,\tau}$. Then \eqref{recurrence-CS} can be written in the following matrix form:
\begin{align}\label{eq-CS-recurrence}
C_{n+1}=\bar{C}_nL_n =
\left(
\begin{array}{cc}
1 & O\\
O & C_n
\end{array}
\right)L_{n},
\end{align}
where
\begin{equation}\label{eq-ln}
L_{n}=\left(
\begin{array}{ccccccc}
1 &\quad &\quad &\quad&\quad&\quad&\quad\\
s_0 & r_0 & \quad &\quad & \quad &\quad&\quad\\
t_1 & s_1 & r_1  &\quad &\quad &\quad&\quad\\
\quad& t_2 & s_2 & r_2 &\quad&\quad&\quad\\
\quad&\quad&\ddots & \ddots &\ddots&\quad&\quad\\
\quad &\quad&\quad&t_{n-1}&s_{n-1}&r_{n-1}&\quad\\
\quad&\quad&\quad&\quad&t_n&s_n&r_n
\end{array}
\right).
\end{equation}
Therefore, in view of \eqref{eq-CS-recurrence} and Lemma \ref{lem-transfer-matrix}, the planar network construction of $C^{\gamma,\sigma,\tau}$ can be reduced to that of $L_n$.

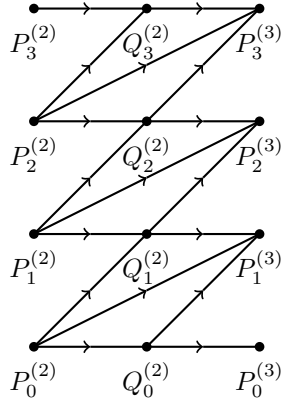
\begin{figure}[ht]
         	\centering
         	\begin{tikzpicture}
         	[place/.style={thick,fill=black!100,circle,inner sep=0pt,minimum size=1mm,draw=black!100},scale=1.5]
        	\draw [thick] [->] (2.5,2.5) -- (3,2.5);
         	\draw [thick] [->] (3,2.5) -- (3.5,2.5) -- (4,2.5);
        	\draw [thick] (4,2.5) -- (4.5,2.5);
         	\node [place,label=below:{\footnotesize$P_3^{(2)}$}] at (2.5,2.5) {};
         	\node [place,label=below:{\footnotesize$Q_3^{(2)}$}] at (3.5,2.5) {};
         	\node [place,label=below:{\footnotesize$P_3^{(3)}$}] at (4.5,2.5) {};
            \draw [thick] [->](2.5,1.5) -- (3,1.5);
         	\draw [thick] [->] (3,1.5) -- (3.5,1.5) -- (4,1.5);
         	\draw [thick] (4,1.5) -- (4.5,1.5);
         	\node [place,label=below:{\footnotesize$P_2^{(2)}$}] at (2.5,1.5) {};
         	\node [place,label=below:{\footnotesize$Q_2^{(2)}$}] at (3.5,1.5) {};
         	\node [place,label=below:{\footnotesize$P_2^{(3)}$}] at (4.5,1.5) {};
         	\draw [thick] [->] (2.5,1.5) --(3,2);
         	\draw [thick] (3,2) --  (3.5,2.5) ;
         	\draw [thick] [->] (2.5,1.5) -- (3.5,2);
         	\draw [thick] (3.5,2) -- (4.5,2.5);
        	\draw [thick] [->] (3.5,1.5) -- (4,2);
        	\draw [thick] (4,2) -- (4.5,2.5);
         	\draw [thick] [->] (2.5,0.5) -- (3,0.5);
         	\draw [thick] [->] (3,0.5) -- (3.5,0.5) -- (4,0.5);
         	\draw [thick] (4,0.5) -- (4.5,0.5);
         	\node [place,label=below:{\footnotesize$P_1^{(2)}$}] at (2.5,0.5) {};
         	\node [place,label=below:{\footnotesize$Q_1^{(2)}$}] at (3.5,0.5) {};
         	\node [place,label=below:{\footnotesize$P_1^{(3)}$}] at (4.5,0.5) {};
         	\draw [thick] [->] (2.5,0.5) -- (3,1);
         	\draw [thick] (3,1) -- (3.5,1.5);
         	\draw [thick] [->] (2.5,0.5) -- (3.5,1);
         	\draw [thick] (3.5,1) -- (4.5,1.5);
         	\draw [thick] [->] (3.5,0.5) -- (4,1);
         	\draw [thick] (4,1) -- (4.5,1.5);
         	\draw [thick] [->] (2.5,-0.5) -- (3,-0.5);
         	\draw [thick] [->] (3,-0.5) -- (3.5,-0.5) -- (4,-0.5);
         	\draw [thick]  (4,-0.5) -- (4.5,-0.5);
         	\node [place,label=below:{\footnotesize$P_0^{(2)}$}] at (2.5,-0.5) {};
         	\node [place,label=below:{\footnotesize$Q_0^{(2)}$}] at (3.5,-0.5) {};
         	\node [place,label=below:{\footnotesize$P_0^{(3)}$}] at (4.5,-0.5) {};
         	\draw [thick] [->] (2.5,-0.5) -- (3,0);
         	\draw [thick] (3,0) -- (3.5,0.5);
         	\draw [thick] [->] (2.5,-0.5) -- (3.5,0);
         	\draw [thick] (3.5,0) -- (4.5,0.5);
         	\draw [thick] [->] (3.5,-0.5) -- (4,0);
         	\draw [thick] (4,0) -- (4.5,0.5);
         	\end{tikzpicture}
         	\caption{Digraph $D^{L_2}$ of $\mathcal{L}_2$}\label{fig-d2l}
         \end{figure}

We proceed to construct the planar network $\mathcal{L}_n$ for $L_n$. We first build the underlying digraph $D^{L_n}$ of $\mathcal{L}_n$. Let $\mathcal{P}^{(n)}=\{P_0^{(n)},\,P_1^{(n)},\,P_2^{(n)},\,\ldots\}$ and $\mathcal{Q}^{(n)}=\{Q_0^{(n)},\,Q_1^{(n)},\,Q_2^{(n)},\,\ldots \}$ be two sets of lattice points in the plane with coordinates
\begin{align*}
P_{i}^{(n)} = (2n,i), \qquad Q_i^{(n)} = (2n+1,i)
\end{align*}
for $i\geq 0$. The vertex set of $D^{L_n}$ is given by
\begin{align*}
\{P_{0}^{(n)},\,P_{1}^{(n)},\,\ldots,\,P_{n+1}^{(n)},\,
 P_{0}^{(n+1)},\,P_{1}^{(n+1)},\,\ldots,\,P_{n+1}^{(n+1)},\,
Q_{0}^{(n)},\,Q_{1}^{(n)},\,\ldots,\,Q_{n+1}^{(n)}\}.
\end{align*}
The arc set of $D^{L_n}$ is composed of the following three kinds of arcs:
\begin{itemize}
\item[(1)] horizontal arcs $P_{k}^{(n)}\to Q_{k}^{(n)}$ and $Q_{k}^{(n)}\to P_{k}^{(n+1)}$ for $0 \le k \le n+1$;
\item[(2)] diagonal arcs $P_{k}^{(n)}\to Q_{k+1}^{(n)}$ and $Q_{k}^{(n)}\to P_{k+1}^{(n+1)}$ for $0 \le k \le n$; and
\item[(3)] super diagonal arcs $P_{k}^{(n)}\to P_{k+1}^{(n+1)}$ for $0 \le k \le n$.
\end{itemize}
It is easy to observe that $D^{L_n}$ is a planar graph, see Figure \ref{fig-d2l} for an illustration.

To construct the planar network $\mathcal{L}_n$, we need to assign a weight to each arc of $D^{L_n}$.
For each case of Theorem \ref{thm-tp-csmatrix}, we will choose a weight function $\mathrm{wt}$ such that each arc $a$ is assigned a $q$-nonnegative weight. For convenience, we set $r_{-1}=t_{0}=0$.
For the first case of Theorem \ref{thm-tp-csmatrix}, let
\begin{align}
\mathrm{wt}(a)=
\begin{cases}
r_{n-k}, & \mbox{ if } a=P_{k}^{(n)}\to Q_{k}^{(n)}
\mbox{ for some } 0\le k\le n
;\\
t_{n-k}, & \mbox{ if } a=P_{k}^{(n)}\to Q_{k+1}^{(n)}
\mbox{ for some } 0\le k\le n
;\\
s_{n-k}-r_{n-k}-t_{n-k}, & \mbox{ if } a= P_{k}^{(n)}\to P_{k+1}^{(n+1)}
\mbox{ for some } 0\le k\le n
;\\
1,& \mbox{ otherwise}.
\end{cases}\label{eq-wt1}
\end{align}
For the second case, let
\begin{align}
\mathrm{wt}(a)=
\begin{cases}
r_{n-k}, & \mbox{ if } a=Q_{k}^{(n)}\to P_{k}^{(n+1)}
\mbox{ for some } 0\le k\le n
;\\
t_{n-k+1}, & \mbox{ if } a=Q_{k}^{(n)}\to P_{k+1}^{(n+1)}
\mbox{ for some } 0\le k\le n
;\\
s_{n-k}-r_{n-k-1}-t_{n-k+1}, & \mbox{ if } a= P_{k}^{(n)}\to P_{k+1}^{(n+1)}
\mbox{ for some } 0\le k\le n
;\\
0,&\mbox{ if } a=P_{n}^{(n)}\to Q_{n+1}^{(n)};\\
1,& \mbox{ otherwise}.
\end{cases}\label{eq-wt2}
\end{align}
For the third case, let
\begin{align}
\mathrm{wt}(a)=
\begin{cases}
r_{n-k}, & \mbox{ if } a=Q_{k}^{(n)}\to P_{k}^{(n+1)}
\mbox{ for some } 0\le k\le n
;\\
t_{n-k}, & \mbox{ if } a=P_{k}^{(n)}\to Q_{k+1}^{(n)}
\mbox{ for some } 0\le k\le n
;\\
s_{n-k}-r_{n-k-1}t_{n-k}-1, & \mbox{ if } a= P_{k}^{(n)}\to P_{k+1}^{(n+1)}
\mbox{ for some } 0\le k\le n
;\\
1,& \mbox{ otherwise}.
\end{cases}\label{eq-wt3}
\end{align}
For the fourth case, let
\begin{align}
\mathrm{wt}(a)=
\begin{cases}
r_{n-k}, & \mbox{ if } a=P_{k}^{(n)}\to Q_{k}^{(n)}\mbox{ for some } 0\le k\le n;\\
t_{n-k+1}, & \mbox{ if } a=Q_{k}^{(n)}\to P_{k+1}^{(n+1)}\mbox{ for some } 0\le k\le n;\\
s_{n-k}-r_{n-k}t_{n-k+1}-1, & \mbox{ if } a= P_{k}^{(n)}\to P_{k+1}^{(n+1)}\mbox{ for some } 0\le k\le n-1;\\
s_0-r_0t_1,&\mbox{ if } a=P_{n}^{(n)}\to P_{n+1}^{(n+1)};\\
0,&\mbox{ if } a=P_{n}^{(n)}\to Q_{n+1}^{(n)};\\
1,& \mbox{ otherwise}.
\end{cases}\label{eq-wt4}
\end{align}
For the fifth case, let
\begin{align}
\mathrm{wt}(a)=
\begin{cases}
b_{n-k}, & \mbox{ if } a=P_{k}^{(n)}\to Q_{k+1}^{(n)}\mbox{ for some } 0\le k\le n;\\
c_{n-k}, & \mbox{ if } a=Q_{k}^{(n)}\to P_{k+1}^{(n+1)}\mbox{ for some } 0\le k\le n;\\
0, & \mbox{ if } a= P_{k}^{(n)}\to P_{k+1}^{(n+1)}\mbox{ for some } 0\le k\le n;\\
1,& \mbox{ otherwise}.
\end{cases}\label{eq-wt5}
\end{align}

Now by the above construction we have the following result.

\begin{prop}\label{prop-network-Ln}
Suppose that $L_n$ is given by \eqref{eq-ln}. Then
$$\mathcal{L}_n = (D^{L_n},\mathrm{wt}_{D^{L_n}},(P_{n+1}^{(n)},\ldots,P_{0}^{(n)}), (P_{n+1}^{(n+1)},\ldots,P_{0}^{(n+1)}))$$
is a planar network for $L_n$, where  $\mathrm{wt}_{D^{L_n}}$ is
given by \eqref{eq-wt1}, \eqref{eq-wt2},
\eqref{eq-wt3}, \eqref{eq-wt4} or \eqref{eq-wt5} if applicable.
Consequently,
\begin{equation*}
L_n=\left(GF_{\mathcal{L}_n}(P_{n+1-i}^{(n)},P_{n+1-j}^{(n+1)})\right)_{0 \le i,j \le n+1}.
\end{equation*}
\end{prop}


Now we are in a position to present the planar network $\mathcal{C}_{n}$ for $C_{n}$, which can be recursively constructed as follows:

\begin{itemize}
	\item For $n=1$, we take $\mathcal{C}_{1}$ to be the planar network $\mathcal{L}_{0}$, whose digraph $D^{C_1}( = D^{L_0})$ is shown as in
Figure \ref{fig-dc1}.
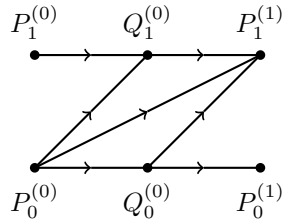
\begin{figure}[h]
\centering
\begin{tikzpicture}
[place/.style={thick,fill=black!100,circle,inner sep=0pt,minimum size=1mm,draw=black!100}, scale = 1.5]
\draw [thick] [->] (-2,2) -- (-1.5,2);
\draw [thick] (-1.5,2) -- (-1,2);
\draw [thick] [->] (-1,2) -- (-0.5,2);
\draw [thick] (-0.5,2) -- (0,2);
\draw [thick] [->] (-1,1)  -- (-0.5,1.5);
\draw [thick] (-0.5,1.5) -- (0,2);
\draw [thick] [->] (-2,1) -- (-1.5,1);
\draw [thick] (-1.5,1) -- (-1,1);
\draw [thick] [->] (-1,1) -- (-0.5,1);
\draw [thick] (-0.5,1) -- (0,1);
\draw [thick] [->] (-2,1) -- (-1.5,1.5);
\draw [thick] (-1.5,1.5) -- (-1,2);
\draw [thick] [->] (-2,1) -- (-1,1.5);
\draw [thick] (-1,1.5) -- (0,2);
\node [place,label=below:{\footnotesize$P_0^{(0)}$}] at (-2,1) {};
\node [place,label=above:{\footnotesize$P_1^{(0)}$}] at (-2,2) {};
\node [place,label=above:{\footnotesize$Q_1^{(0)}$}] at (-1,2) {};
\node [place,label=above:{\footnotesize$P_1^{(1)}$}] at (0,2) {};
\node [place,label=below:{\footnotesize$Q_0^{(0)}$}] at (-1,1) {};
\node [place,label=below:{\footnotesize$P_0^{(1)}$}] at (0,1) {};
\end{tikzpicture}
\caption{Digraph $D^{C_1}$ of $\mathcal{C}_1$}\label{fig-dc1}
\end{figure}

\item Provided that $\mathcal{C}_{n}$
has been constructed for some $n\geq 1$, we continue to build $\mathcal{C}_{n+1}$.
Let $V(D^{C_n})$ and $A(D^{C_n})$ denote the vertex set and arc set of $D^{C_n}$, respectively. By adding the vertices $P_{n+1}^{(0)},P_{n+1}^{(1)},\ldots,P_{n+1}^{(n)}$ to $V(D^{C_n})$ and arcs $P_{n+1}^{(i)} \to P_{n+1}^{(i+1)}$ ($0\leq i\leq n-1$)  to $A(D^{C_n})$, we obtain $D^{\bar{C}_n}$, the digraph of $\bar{\mathcal{C}}_n$. Let $\mathrm{wt}_{D^{\bar{C}_n}}(a)$ be equal to $\mathrm{wt}_{D^{C_n}}(a)$ for $a \in A(D^{C_n})$ and equal to 1 otherwise. Then
\[
    \mathcal{\bar{C}}_n = (D^{\bar{C}_n},\mathrm{wt}_{D^{\bar{C}_n}}, (P_{n+1}^{(0)},P_n^{(0)},\ldots,P_0^{(0)}), (P_{n+1}^{(n)},P_n^{(n)},\ldots,P_0^{(n)}))
\]
is a network for $\bar{C}_n$.
In view of \eqref{eq-CS-recurrence}, the planar network $\mathcal{C}_{n+1}$ for $C_{n+1}$ can be obtained from $\mathcal{\bar{C}}_n$ and $\mathcal{L}_n$ by Lemma \ref{lem-transfer-matrix}. The digraph $D^{C_3}$ is depicted in Figure \ref{fig-d3c}.

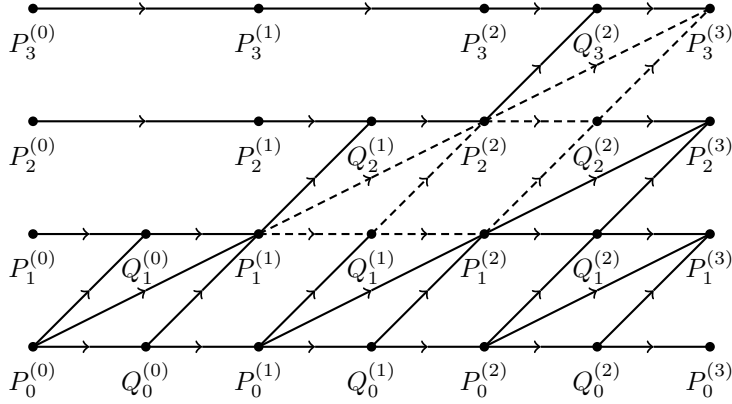
\begin{figure}[ht]
\centering
\begin{tikzpicture}
[place/.style={thick,fill=black!100,circle,inner sep=0pt,minimum size=1mm,draw=black!100},scale=1.5]
\draw [thick] [->] (-1.5,2.5) -- (-0.5,2.5);
\draw [thick] [->] (-0.5,2.5) -- (1.5,2.5);
\draw [thick] (1.5,2.5) -- (2.5,2.5);
\draw [thick] [->] (2.5,2.5) -- (3,2.5);
\draw [thick] (3,2.5) -- (3.5,2.5);
\draw [thick] [->] (3.5,2.5) -- (4,2.5);
\draw [thick] (4,2.5) -- (4.5,2.5);
\node [place,label=below:{\footnotesize$P_3^{(0)}$}] at (-1.5,2.5) {};
\node [place,label=below:{\footnotesize$P_3^{(1)}$}] at (0.5,2.5) {};
\node [place,label=below:{\footnotesize$P_3^{(2)}$}] at (2.5,2.5) {};
\node [place,label=below:{\footnotesize$Q_3^{(2)}$}] at (3.5,2.5) {};
\node [place,label=below:{\footnotesize$P_3^{(3)}$}] at (4.5,2.5) {};
\draw [thick] [->] (-1.5,1.5) -- (-0.5,1.5);
\draw [thick] [->] (-0.5,1.5) -- (1,1.5);
\draw [thick] [->] (1,1.5) -- (2,1.5);
\draw [thick] (2,1.5) --(2.5,1.5);
\draw [thick][->][densely dashed] (2.5,1.5) -- (3,1.5);
\draw [thick] [densely dashed](3,1.5) --(3.5,1.5);
\draw [thick][->] (3.5,1.5) -- (4,1.5);
\draw [thick] (4,1.5) -- (4.5,1.5);
\node [place,label=below:{\footnotesize$P_2^{(0)}$}] at (-1.5,1.5) {};
\node [place,label=below:{\footnotesize$P_2^{(1)}$}] at (0.5,1.5) {};
\node [place,label=below:{\footnotesize$Q_2^{(1)}$}] at (1.5,1.5) {};
\node [place,label=below:{\footnotesize$P_2^{(2)}$}] at (2.5,1.5) {};
\node [place,label=below:{\footnotesize$Q_2^{(2)}$}] at (3.5,1.5) {};
\node [place,label=below:{\footnotesize$P_2^{(3)}$}] at (4.5,1.5) {};
\draw [thick] [->] (2.5,1.5) --(3,2);
\draw [thick] (3,2) --  (3.5,2.5) ;
\draw [thick] [->][densely dashed] (2.5,1.5) -- (3.5,2);
\draw [thick] [densely dashed] (3.5,2) -- (4.5,2.5);
\draw [thick] [->] [densely dashed](3.5,1.5) -- (4,2);
\draw [thick][densely dashed] (4,2) -- (4.5,2.5);
\draw [thick] [->] (-1.5,0.5) -- (-1,0.5);
\draw [thick] [->] (-1,0.5) --  (0,0.5);
\draw [thick] (0,0.5) -- (0.5,0.5);
\draw [thick] [->][densely dashed] (0.5,0.5) -- (1,0.5);
\draw [thick] [->][densely dashed] (1,0.5) -- (2,0.5);
\draw [thick] [densely dashed] (2,0.5) -- (2.5,0.5);
\draw [thick] [->] (2.5,0.5) -- (3,0.5);
\draw [thick] [->] (3,0.5) -- (4,0.5);
\draw [thick] (4,0.5) -- (4.5,0.5);
\node [place,label=below:{\footnotesize$P_1^{(0)}$}] at (-1.5,0.5) {};
\node [place,label=below:{\footnotesize$Q_1^{(0)}$}] at (-0.5,0.5) {};
\node [place,label=below:{\footnotesize$P_1^{(1)}$}] at (0.5,0.5) {};
\node [place,label=below:{\footnotesize$Q_1^{(1)}$}] at (1.5,0.5) {};
\node [place,label=below:{\footnotesize$P_1^{(2)}$}] at (2.5,0.5) {};
\node [place,label=below:{\footnotesize$Q_1^{(2)}$}] at (3.5,0.5) {};
\node [place,label=below:{\footnotesize$P_1^{(3)}$}] at (4.5,0.5) {};
\draw [thick] [->] (0.5,0.5) --(1,1);
\draw [thick] (1,1) -- (1.5,1.5);
\draw [thick][densely dashed] (1.5,1) -- (2.5,1.5);
\draw [thick] [->] [densely dashed] (1.5,0.5) -- (2,1);
\draw [thick] [densely dashed] (2,1) -- (2.5,1.5);
\draw [thick] [->] [densely dashed] (2.5,0.5) -- (3,1);
\draw [thick] [densely dashed] (3,1) -- (3.5,1.5);
\draw [thick] [->] (2.5,0.5) -- (3.5,1);
\draw [thick] (3.5,1) -- (4.5,1.5);
\draw [thick] [->] (3.5,0.5) -- (4,1);
\draw [thick] (4,1) -- (4.5,1.5);
\draw [thick] [->] (-1.5,-0.5) -- (-1,-0.5);
\draw [thick] [->] (-1,-0.5) -- (-0.5,-0.5) -- (0,-0.5);
\draw [thick] [->] (0,-0.5) -- (1,-0.5);
\draw [thick] [->] (1,-0.5) -- (1.5,-0.5) -- (2,-0.5);
\draw [thick] [->] (2,-0.5) -- (2.5,-0.5) -- (3,-0.5);
\draw [thick] [->] (3,-0.5) -- (3.5,-0.5) -- (4,-0.5);
\draw [thick] (4,-0.5) -- (4.5,-0.5);
\node [place,label=below:{\footnotesize$P_0^{(0)}$}] at (-1.5,-0.5) {};
\node [place,label=below:{\footnotesize$Q_0^{(0)}$}] at (-0.5,-0.5) {};
\node [place,label=below:{\footnotesize$P_0^{(1)}$}] at (0.5,-0.5) {};
\node [place,label=below:{\footnotesize$Q_0^{(1)}$}] at (1.5,-0.5) {};
\node [place,label=below:{\footnotesize$P_0^{(2)}$}] at (2.5,-0.5) {};
\node [place,label=below:{\footnotesize$Q_0^{(2)}$}] at (3.5,-0.5) {};
\node [place,label=below:{\footnotesize$P_0^{(3)}$}] at (4.5,-0.5) {};
\draw [thick] [->] (-0.5,-0.5)  -- (0,0);
\draw [thick] (0,0) --  (0.5,0.5);
\draw [thick] [->] (0.5,-0.5) -- (1,0);
\draw [thick] (1,0) --(1.5,0.5);
\draw [thick] [->] (0.5,-0.5) -- (1.5,0);
\draw [thick] (1.5,0) -- (2.5,0.5);
\draw [thick] [->]  (1.5,-0.5) -- (2,0);
\draw [thick] (2,0) -- (2.5,0.5);
\draw [thick] [->]  (2.5,-0.5) -- (3,0);
\draw [thick] (3,0) -- (3.5,0.5);
\draw [thick] [->] (2.5,-0.5) -- (3.5,0);
\draw [thick] (3.5,0) -- (4.5,0.5);
\draw [thick] [->] (3.5,-0.5) -- (4,0);
\draw [thick] (4,0) -- (4.5,0.5);
\draw [thick] [->] (-1.5,-0.5)  -- (-1,0);
\draw [thick] (-1,0) -- (-0.5,0.5);
\draw [thick] [->] (-1.5,-0.5) -- (-0.5,0);
\draw [thick]  (-0.5,0) -- (0.5,0.5);
\draw [thick][densely dashed][->] (0.5,0.5) -- (1.5,1);
\end{tikzpicture}
\caption{Digraph $D^{C_3}$ of $\mathcal{C}_3$ and an example of $D^{H_1}$}\label{fig-d3c}
\end{figure}
\end{itemize}

Note that the above construction implies that the weight function $\mathrm{wt}_{D^{C_n}}$ of $\mathcal{C}_n$ is uniquely determined by the weight functions of  
$\mathcal{L}_0,\mathcal{L}_1,\ldots,\mathcal{L}_{n-1}$.
Therefore, we may use $(w_0,w_1,\ldots,w_{n-1})$ to represent the weight function $\mathrm{wt}_{D^{C_n}}$, where $w_i$ denotes the weight function of $\mathcal{L}_i$. Applying Lemma \ref{lem-transfer-matrix}, Proposition \ref{prop-network-Ln} and the relation \eqref{eq-CS-recurrence}, we have the following result.

\begin{thm}\label{thm-network-CS}
Let $C_n$ be the $n$th leading principal submatrix of $C^{\gamma,\sigma,\tau}$.
Then
\begin{equation}\label{eq-network-for-cn}
\mathcal{C}_n = (D^{C_n},\mathrm{wt}_{D^{C_n}},(P_{n}^{(0)},\ldots,P_{0}^{(0)}), (P_{n}^{(n)},\ldots,P_{0}^{(n)}))
\end{equation}
is a planar network for $C_n$, where
$\mathrm{wt}_{D^{C_n}}=(w_0,\,w_1,\,\ldots,\,w_{n-1})$
and each $w_i$ is given by one of \eqref{eq-wt1}, \eqref{eq-wt2},
\eqref{eq-wt3}, \eqref{eq-wt4} or \eqref{eq-wt5} if applicable.
Consequently, 
\begin{equation}\label{eq-cn-gf}
    C_{n}= \left(GF_{\mathcal{C}_n}(P_{n-i}^{(0)},P_{n-j}^{(n)})\right)_{0 \le i,j \le n}.
\end{equation}
\end{thm}


\noindent \textbf{Remark.} Our construction of $C^{\gamma,\sigma,\tau}$ is different from that of Pan and Zeng \cite{Pan-Zeng-2016}. In their construction, there is exactly one type of weight functions for each planar network. However, in our network $\mathcal{C}_n$, we may choose different types of weight functions for $\mathcal{L}_0,\mathcal{L}_1,\ldots,\mathcal{L}_{n-1}$.
Actually, if we choose the same type of weight functions for all $w_i$ in $\mathrm{wt}_{D^{C_n}}$, then Pan and Zeng's network for $C_n$ is isomorphic to
\begin{equation*}
(D^{C'_n},\mathrm{wt}_{D^{C'_n}},(P_{n}^{(n)},P_{n-1}^{(n-1)},\ldots,P_{0}^{(0)}), (P_{n}^{(n)},P_{n-1}^{(n)},\ldots,P_{0}^{(n)})),
\end{equation*}
where $D^{C'_n}$ is obtained from  $D^{C_n}$ by removing all the vertices and arcs strictly above the directed path $P_{0}^{(0)}\to P_{1}^{(1)}\to \cdots\to P_{n}^{(n)}$, and $\mathrm{wt}_{D^{C'_n}}$
is the restriction of $\mathrm{wt}_{D^{C_n}}$ to $D^{C'_n}$.

Given the planar network interpretation for Catalan-Stieltjes matrix $C^{\gamma,\sigma,\tau}$, we have the following combinatorial interpretation for the corresponding Hankel matrix $H^{\gamma,\sigma,\tau}$.

\begin{cor}\label{thm-network-CS-Hankel}
Fix $k \ge 0$. Let $H_n$ be the $n$th leading principal submatrix of $H^{\gamma,\sigma,\tau}$. Let $D^{H_n}$ be the subgraph of $D^{C_{2n+k}}$ induced by the union of arcs of all directed paths from $P_k^{(k)}$ to $P_{2n+k}^{(2n+k)}$, and let $\mathrm{wt}_{D^{H_n}}$ be the restriction of $\mathrm{wt}_{D^{C_{2n+k}}}$ to $D^{H_n}$. Then
\begin{equation*}
\mathcal{H}_n = (D^{H_n},\mathrm{wt}_{D^{H_n}},(P_{n+k}^{(n+k)},\ldots,P_{k}^{(k)}), (P_{n+k}^{(n+k)},\ldots,P_{2n+k}^{(2n+k)}))
\end{equation*}
is a planar network for $H_n$. Consequently,
\begin{equation*}
H_n = \left(GF_{\mathcal{H}_n}(P_{n+k-i}^{(n+k-i)},P_{n+k+j}^{(n+k+j)})\right)_{0 \le i,j \le n},
\end{equation*}
\end{cor}

\pf It suffices to prove that $$GF_{\mathcal{H}_n}(P_{n+k-i}^{(n+k-i)}, P_{n+k+j}^{(n+k+j)}) = c_{i+j,0}$$
for any $n,k \ge 0$ and $0 \le i,j \le n$.
Substituting $n+k+j$ for $n$ in  \eqref{eq-cn-gf} leads to
$$c_{i+j,0} =GF_{\mathcal{C}_{n+k+j}}(P_{n+k-i}^{(0)}, P_{n+k+j}^{(n+k+j)}).$$
Note that there is exactly one directed path from $P_{n+k-i}^{(0)} \to P_{n+k-i}^{(n+k-i)}$ and the weights of all the arcs in this path are 1.
Thus we have   $$GF_{\mathcal{C}_{n+k+j}}(P_{n+k-i}^{(0)}, P_{n+k+j}^{(n+k+j)}) = GF_{\mathcal{C}_{n+k+j}}(P_{n+k-i}^{(n+k-i)}, P_{n+k+j}^{(n+k+j)}).$$
 Observe that $D^{H_n}$ and $D^{C_{n+k+j}}$ are both subgraph of $D^{C_{2n+k}}$.
Moreover, $\mathrm{wt}_{D^{H_n}}$ and $\mathrm{wt}_{D^{C_{n+k+j}}}$ can be considered as restrictions of $\mathrm{wt}_{D^{C_{2n+k}}}$ to $D^{H_n}$ and $D^{C_{n+k+j}}$, respectively, and hence
\[
GF_{\mathcal{C}_{n+k+j}}(P_{n+k-i}^{(n+k-i)}, P_{n+k+j}^{(n+k+j)}) = GF_{\mathcal{H}_{n}}(P_{n+k-i}^{(n+k-i)}, P_{n+k+j}^{(n+k+j)}).
\]
This completes the proof.
\qed


Taking $k=1$ and $n=1$, the digraph $D^{H_1}$ in Corollary \ref{thm-network-CS-Hankel} is shown as the dashed part of Figure \ref{fig-d3c}.

When $r_k = 1$ for $k \ge 0$, we can give another planar network interpretation for the Hankel matrix. In this case, we denote by $\tilde{H}_n$ the $n$th leading principal submatrix of $H^{\gamma,\sigma,\tau}$. It is known that  \begin{equation}
\tilde{H}_n=C_nT_nC_n^{T},
\end{equation}
where $C_n^T$ denotes the transpose of $C_n$ and
\begin{equation*}
T_n=\left(
\begin{array}{ccccc}
1&\quad&\quad&\quad&\quad\\
\quad&t_1&\quad&\quad&\quad\\
\quad&\quad&t_1t_2&\quad&\quad\\
\quad&\quad&\quad&\ddots&\quad\\
\quad&\quad&\quad&\quad&t_1t_2\cdots t_n
\end{array}
\right),
\end{equation*}
see \cite{Aigner-2001} for more details.

Given the planar network $\mathcal{C}_n$,
it is straightforward to generate a planar network for $C_n^T$. Specifically, let $D^{C_n^T}$ be the digraph obtained from $D^{C_n}$ by taking its mirror image with respect to the vertical line $x=2n+1/2$ in the plane and then reversing all the arcs. The weight function $\mathrm{wt}_{D^{C_n^T}}$ of $D^{C_n^T}$
assigns to each arc the weight of its preimage in $\mathcal{C}_n$. Denote the mirror image of $P_{i}^{(j)}$ (respectively, $Q_{i}^{(j)}$) by $\bar{P}_{i}^{(j)}$ (respectively, $\bar{Q}_{i}^{(j)}$).
Now
\begin{equation}\label{eq-network-for-cnt}
\mathcal{C}^T_n = (D^{C_n^T},\mathrm{wt}_{D^{C_n^T}}, (\bar{P}_{n}^{(n)},\ldots,\bar{P}_{0}^{(n)}),
(\bar{P}_{n}^{(0)},\ldots,\bar{P}_{0}^{(0)}))
\end{equation}
is a planar network for $C_n^T$.
Let $D^{T_n}$ denote the digraph with
vertex set
$$\{P_{0}^{(n)} ,\ldots,P_{n}^{(n)},\bar{P}_{0}^{(n)} ,\ldots,\bar{P}_{n}^{(n)}\}$$ and arc set
$$\{P_{0}^{(n)}\to \bar{P}_{0}^{(n)},\ldots,
P_{n}^{(n)}\to \bar{P}_{n}^{(n)}\}.$$
Let $\mathrm{wt}_{D^{T_n}}$ be the weight function which maps each arc $P_{i}^{(n)} \to \bar{P}_{i}^{(n)}$ to $t_1\cdots t_{n-i}$ for any $0\leq i\leq n$. Then
\begin{equation}\label{eq-network-for-tn}
\mathcal{T}_n = (D^{T_n},\mathrm{wt}_{D^{T_n}}, ({P}_{n}^{(n)},\ldots,{P}_{0}^{(n)}),
(\bar{P}_{n}^{(n)},\ldots,\bar{P}_{0}^{(n)}))
\end{equation}
is a planar network for $T_n$.
We immediately obtain the following result.

\begin{cor}
Suppose that $\mathcal{C}_n$, $\mathcal{T}_n$ and $\mathcal{C}^T_n$ are respectively given by \eqref{eq-network-for-cn},
\eqref{eq-network-for-tn} and \eqref{eq-network-for-cnt}.
Let
\begin{equation*}
\mathcal{\tilde{H}}_n =
(D^{C_nT_nC_n^T},\mathrm{wt}_{D^{C_nT_nC_n^T}},(P_{n}^{(0)},\ldots,P_{0}^{(0)}), (\bar{P}_{n}^{(0)},\ldots,\bar{P}_{0}^{(0)}))
\end{equation*}
be the planar network as constructed in  Lemma \ref{lem-transfer-matrix} from  $\mathcal{C}_n$, $\mathcal{T}_n$ and $\mathcal{C}^T_n$. Then $\mathcal{\tilde{H}}_n$
is a planar network for $\tilde{H}_n$. Consequently,
\begin{equation*}
\tilde{H}_n = \left(GF_{\mathcal{\tilde{H}}_n}(P_{n-i}^{(0)},\bar{P}_{n-j}^{(0)})\right)_{0 \le i,j \le n}.
\end{equation*}
\end{cor}

\begin{figure}[htb]
	\centering
\begin{tikzpicture}
	[place/.style={thick,fill=black!100,circle,inner sep=0pt,minimum size=1mm,draw=black!100},scale=1.15]
	\draw [thick] [->] (-1.5,2.5) -- (-0.5,2.5);
    \draw [thick]  (-0.5,2.5) -- (0.5,2.5);
	\draw [thick] [->] (0.5,2.5) -- (1.5,2.5);
    \draw [thick]  (1.5,2.5) -- (2.5,2.5);
	\draw [thick] [->]  (2.5,2.5) -- (3,2.5);
	\draw [thick] [->] (3,2.5) -- (3.5,2.5) -- (4,2.5);
	\draw [thick] [->] (4,2.5) -- (4.5,2.5) -- (5,2.5);
	\draw [thick] [->] (5,2.5) -- (5.5,2.5) -- (6,2.5);
	\draw [thick] [->] (6,2.5) -- (6.5,2.5) -- (7,2.5);
	\draw [thick] [->]  (7,2.5) -- (8.5,2.5);
    \draw [thick] (8.5,2.5) -- (9.5,2.5);
	\draw [thick] [->] (9.5,2.5) -- (10.5,2.5);
    \draw [thick] (10.5,2.5) -- (11.5,2.5);
	\node [place,label=below:{\tiny$P_3^{(0)}$}] at (-1.5,2.5) {};
	\node [place,label=below:{\tiny$P_3^{(1)}$}] at (0.5,2.5) {};
	\node [place,label=below:{\tiny$P_3^{(2)}$}] at (2.5,2.5) {};
	\node [place,label=below:{\tiny$Q_3^{(2)}$}] at (3.5,2.5) {};
	\node [place,label=below:{\tiny$P_3^{(3)}$}] at (4.5,2.5) {};
	\node [place,label=below:{\tiny$\bar{P}_3^{(3)}$}] at (5.5,2.5) {};
	\node [place,label=below:{\tiny$\bar{Q}_3^{(2)}$}] at (6.5,2.5) {};
	\node [place,label=below:{\tiny$\bar{P}_3^{(2)}$}] at (7.5,2.5) {};
	\node [place,label=below:{\tiny$\bar{P}_3^{(1)}$}] at (9.5,2.5) {};
	\node [place,label=below:{\tiny$\bar{P}_3^{(0)}$}] at (11.5,2.5) {};
	\draw [thick] [->] (5.5,2.5) --(6,2);
	\draw [thick] (6,2) -- (6.5,1.5);
	\draw [thick] [->] (5.5,2.5) -- (6.5,2);
	\draw [thick] (6.5,2) -- (7.5,1.5);
	\draw [thick] [->] (6.5,2.5) -- (7,2);
	\draw [thick] (7,2) -- (7.5,1.5);
	\draw [thick] [->] (-1.5,1.5) -- (-0.5,1.5);
	\draw [thick]  (-0.5,1.5) -- (0.5,1.5);
	\draw [thick] [->]  (0.5,1.5) -- (1,1.5);
	\draw [thick] [->] (1,1.5) -- (1.5,1.5) -- (2,1.5);
	\draw [thick] [->] (2,1.5) -- (2.5,1.5) -- (3,1.5);
	\draw [thick] [->] (3,1.5) -- (3.5,1.5) -- (4,1.5);
	\draw [thick] [->] (4,1.5) -- (4.5,1.5) -- (5,1.5);
	\draw [thick] [->] (5,1.5) -- (5.5,1.5) -- (6,1.5);
	\draw [thick] [->] (6,1.5) -- (6.5,1.5) -- (7,1.5);
	\draw [thick] [->] (7,1.5) -- (7.5,1.5) -- (8,1.5);
	\draw [thick] [->] (8,1.5) -- (8.5,1.5) -- (9,1.5);
	\draw [thick] [->] (9,1.5) --  (10.5,1.5);
	\draw [thick] (10.5,1.5) --  (11.5,1.5);
	\node [place,label=below:{\tiny$P_2^{(0)}$}] at (-1.5,1.5) {};
	\node [place,label=below:{\tiny$P_2^{(1)}$}] at (0.5,1.5) {};
	\node [place,label=below:{\tiny$Q_2^{(1)}$}] at (1.5,1.5) {};
	\node [place,label=below:{\tiny$P_2^{(2)}$}] at (2.5,1.5) {};
	\node [place,label=below:{\tiny$Q_2^{(2)}$}] at (3.5,1.5) {};
	\node [place,label=below:{\tiny$P_2^{(3)}$}] at (4.5,1.5) {};
	\node [place,label=below:{\tiny$\bar{P}_2^{(3)}$}] at (5.5,1.5) {};
	\node [place,label=below:{\tiny$\bar{Q}_2^{(2)}$}] at (6.5,1.5) {};
	\node [place,label=below:{\tiny$\bar{P}_2^{(2)}$}] at (7.5,1.5) {};
	\node [place,label=below:{\tiny$\bar{Q}_2^{(1)}$}] at (8.5,1.5) {};
	\node [place,label=below:{\tiny$\bar{P}_2^{(1)}$}] at (9.5,1.5) {};
	\node [place,label=below:{\tiny$\bar{P}_2^{(0)}$}] at (11.5,1.5) {};
	\draw [thick] [->] (2.5,1.5) --(3,2);
	\draw [thick]  (3,2) --  (3.5,2.5) ;
	\draw [thick] [->] (2.5,1.5) -- (3.5,2);
	\draw [thick] (3.5,2) -- (4.5,2.5);
	\draw [thick] [->] (3.5,1.5) -- (4,2);
	\draw [thick] (4,2) -- (4.5,2.5);
	\draw [thick] [->] (5.5,1.5) --(6,1);
	\draw [thick] (6,1) -- (6.5,0.5);
	\draw [thick] [->] (5.5,1.5) -- (6.5,1);
	\draw [thick] (6.5,1) -- (7.5,0.5);
	\draw [thick] [->] (6.5,1.5) -- (7,1);
	\draw [thick] (7,1) -- (7.5,0.5);
	\draw [thick] [->] (7.5,1.5) --(8,1);
	\draw [thick] (8,1) -- (8.5,0.5);
	\draw [thick] [->] (7.5,1.5) -- (8.5,1);
	\draw [thick] (8.5,1) -- (9.5,0.5);
	\draw [thick] [->] (8.5,1.5) -- (9,1);
	\draw [thick] (9,1) -- (9.5,0.5);
	\draw [thick] [->] (-1.5,0.5) -- (-1,0.5);
	\draw [thick] [->] (-1,0.5) -- (-0.5,0.5) -- (0,0.5);
	\draw [thick] [->] (0,0.5) -- (0.5,0.5) -- (1,0.5);
	\draw [thick] [->] (1,0.5) -- (1.5,0.5) -- (2,0.5);
	\draw [thick] [->] (2,0.5) -- (2.5,0.5) -- (3,0.5);
	\draw [thick] [->] (3,0.5) -- (3.5,0.5) -- (4,0.5);
	\draw [thick] [->] (4,0.5) -- (4.5,0.5) -- (5,0.5);
	\draw [thick] [->] (5,0.5) -- (5.5,0.5) -- (6,0.5);
	\draw [thick] [->] (6,0.5) -- (6.5,0.5) -- (7,0.5);
	\draw [thick] [->] (7,0.5) -- (7.5,0.5) -- (8,0.5);
	\draw [thick] [->] (8,0.5) -- (8.5,0.5) -- (9,0.5);
	\draw [thick] [->] (9,0.5) -- (9.5,0.5) -- (10,0.5);
	\draw [thick] (10,0.5) -- (10.5,0.5);
	\draw [thick] [->] (10.5,0.5) -- (11,0.5);
	\draw [thick] (11,0.5) -- (11.5,0.5);
	\node [place,label=below:{\tiny$P_1^{(0)}$}] at (-1.5,0.5) {};
	\node [place,label=below:{\tiny$Q_1^{(0)}$}] at (-0.5,0.5) {};
	\node [place,label=below:{\tiny$P_1^{(1)}$}] at (0.5,0.5) {};
	\node [place,label=below:{\tiny$Q_1^{(1)}$}] at (1.5,0.5) {};
	\node [place,label=below:{\tiny$P_1^{(2)}$}] at (2.5,0.5) {};
	\node [place,label=below:{\tiny$Q_1^{(2)}$}] at (3.5,0.5) {};
	\node [place,label=below:{\tiny$P_1^{(3)}$}] at (4.5,0.5) {};
	\node [place,label=below:{\tiny$\bar{P}_1^{(3)}$}] at (5.5,0.5) {};
	\node [place,label=below:{\tiny$\bar{Q}_1^{(2)}$}] at (6.5,0.5) {};
	\node [place,label=below:{\tiny$\bar{P}_1^{(2)}$}] at (7.5,0.5) {};
	\node [place,label=below:{\tiny$\bar{Q}_1^{(1)}$}] at (8.5,0.5) {};
	\node [place,label=below:{\tiny$\bar{P}_1^{(1)}$}] at (9.5,0.5) {};
	\node [place,label=below:{\tiny$\bar{Q}_1^{(0)}$}] at (10.5,0.5) {};
	\node [place,label=below:{\tiny$\bar{P}_1^{(0)}$}] at (11.5,0.5) {};
	\draw [thick] [->] (0.5,0.5) --(1,1);
	\draw [thick] (1,1) -- (1.5,1.5);
	\draw [thick] [->] (0.5,0.5) -- (1.5,1);
	\draw [thick] (1.5,1) -- (2.5,1.5);
	\draw [thick] [->] (1.5,0.5) -- (2,1);
	\draw [thick] (2,1) -- (2.5,1.5);
	\draw [thick] [->] (2.5,0.5) -- (3,1);
	\draw [thick] (3,1) -- (3.5,1.5);
	\draw [thick] [->] (2.5,0.5) -- (3.5,1);
	\draw [thick] (3.5,1) -- (4.5,1.5);
	\draw [thick] [->] (3.5,0.5) -- (4,1);
	\draw [thick] (4,1) -- (4.5,1.5);
	\draw [thick] [->] (-1.5,-0.5) -- (-1,-0.5);
	\draw [thick] [->] (-1,-0.5) --  (0,-0.5);
	\draw [thick] [->] (0,-0.5) -- (1,-0.5);
	\draw [thick] [->] (1,-0.5) -- (2,-0.5);
	\draw [thick] [->] (2,-0.5) --  (3,-0.5);
	\draw [thick] [->] (3,-0.5) --  (4,-0.5);
	\draw [thick] [->] (4,-0.5) --  (5,-0.5);
	\draw [thick] [->] (5,-0.5) --  (6,-0.5);
	\draw [thick] [->] (6,-0.5) -- (7,-0.5);
	\draw [thick] [->] (7,-0.5) --  (8,-0.5);
	\draw [thick] [->] (8,-0.5) --  (9,-0.5);
	\draw [thick] [->] (9,-0.5) -- (10,-0.5);
	\draw [thick] (10,-0.5) -- (10.5,-0.5);
	\draw [thick] [->] (10.5,-0.5) -- (11,-0.5);
	\draw [thick] (11,-0.5) -- (11.5,-0.5);
	\node [place,label=below:{\tiny$Q_0^{(0)}$}] at (-0.5,-0.5) {};
	\node [place,label=below:{\tiny$P_0^{(0)}$}] at (-1.5,-0.5) {};
	\node [place,label=below:{\tiny$P_0^{(1)}$}] at (0.5,-0.5) {};
	\node [place,label=below:{\tiny$Q_0^{(1)}$}] at (1.5,-0.5) {};
	\node [place,label=below:{\tiny$P_0^{(2)}$}] at (2.5,-0.5) {};
	\node [place,label=below:{\tiny$Q_0^{(2)}$}] at (3.5,-0.5) {};
	\node [place,label=below:{\tiny$P_0^{(3)}$}] at (4.5,-0.5) {};
	\node [place,label=below:{\tiny$\bar{P}_0^{(3)}$}] at (5.5,-0.5) {};
	\node [place,label=below:{\tiny$\bar{Q}_0^{(2)}$}] at (6.5,-0.5) {};
	\node [place,label=below:{\tiny$\bar{P}_0^{(2)}$}] at (7.5,-0.5) {};
	\node [place,label=below:{\tiny$\bar{Q}_0^{(1)}$}] at (8.5,-0.5) {};
	\node [place,label=below:{\tiny$\bar{P}_0^{(1)}$}] at (9.5,-0.5) {};
	\node [place,label=below:{\tiny$\bar{Q}_0^{(0)}$}] at (10.5,-0.5) {};
	\node [place,label=below:{\tiny$\bar{P}_0^{(0)}$}] at (11.5,-0.5) {};
	\draw [thick] [->] (-1.5,-0.5)  -- (-1,0);
	\draw [thick] (-1,0) --  (-0.5,0.5);
	\draw [thick] [->] (-1.5,-0.5)  -- (-0.5,0);
	\draw [thick] (-0.5,0) --  (0.5,0.5);
	\draw [thick] [->] (-0.5,-0.5)  -- (0,0);
	\draw [thick] (0,0) --  (0.5,0.5);
	\draw [thick] [->] (0.5,-0.5) -- (1,0);
	\draw [thick] (1,0) --(1.5,0.5);
	\draw [thick] [->] (0.5,-0.5) -- (1.5,0);
	\draw [thick] (1.5,0) -- (2.5,0.5);
	\draw [thick] [->]  (1.5,-0.5) -- (2,0);
	\draw [thick] (2,0) -- (2.5,0.5);
	\draw [thick] [->] (2.5,-0.5) -- (3,0);
	\draw [thick] (3,0) -- (3.5,0.5);
	\draw [thick] [->] (2.5,-0.5) -- (3.5,0);
	\draw [thick] (3.5,0) -- (4.5,0.5);
	\draw [thick] [->] (3.5,-0.5) -- (4,0);
	\draw [thick] (4,0) -- (4.5,0.5);
	
	\draw [thick] [->] (5.5,0.5)  -- (6,0);
	\draw [thick] (6,0) --  (6.5,-0.5);
	\draw [thick] [->] (5.5,0.5) -- (6.5,0);
	\draw [thick] (6.5,0) --(7.5,-0.5);
	\draw [thick] [->] (6.5,0.5)  -- (7,0);
	\draw [thick] (7,0) --  (7.5,-0.5);
	\draw [thick] [->] (7.5,0.5)  -- (8,0);
	\draw [thick] (8,0) --  (8.5,-0.5);
	\draw [thick] [->] (7.5,0.5) -- (8.5,0);
	\draw [thick] (8.5,0) --(9.5,-0.5);
	\draw [thick] [->] (8.5,0.5)  -- (9,0);
	\draw [thick] (9,0) --  (9.5,-0.5);
	\draw [thick] [->] (9.5,0.5)  -- (10,0);
	\draw [thick] (10,0) --  (10.5,-0.5);
	\draw [thick] [->] (9.5,0.5)  -- (10.5,0);
	\draw [thick] (10.5,0) --  (11.5,-0.5);
	\draw [thick] [->] (10.5,0.5)  -- (11,0);
	\draw [thick] (11,0) --  (11.5,-0.5);
\end{tikzpicture}
	\caption{Digraph $D^{\tilde{H}_3}$ of $\mathcal{\tilde{H}}_3$}\label{fig-d3h}
\end{figure}
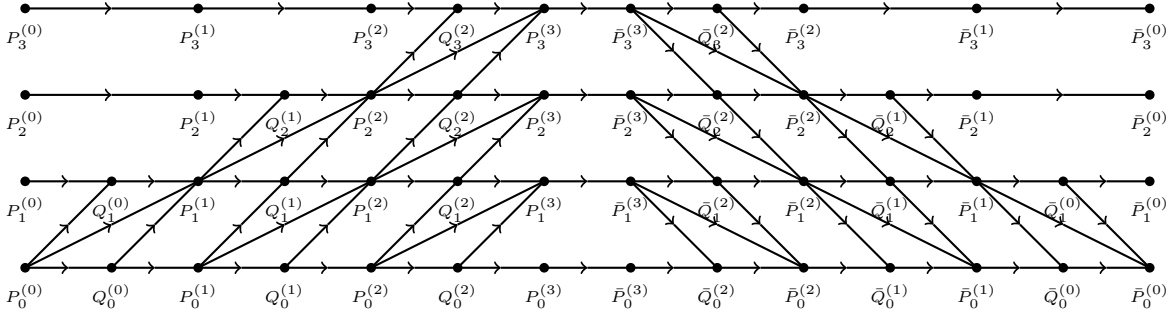

Figure \ref{fig-d3h} presents the digraph $D^{\tilde{H}_3}$ of $\mathcal{\tilde{H}}_3$.

\section{Immanant positivity for  \texorpdfstring{$C^{\gamma,\sigma,\tau}$}{} and \texorpdfstring{$H^{\gamma,\sigma,\tau}$}{}}\label{sect-im}

The main objective of this section is to prove
Theorem \ref{thm-main}. Our proof is based on the planar network interpretation for $C^{\gamma,\sigma,\tau}$ and $H^{\gamma,\sigma,\tau}$ given in Section \ref{sect-ci}. At the end of this section, we will also present a stronger result than
Theorem \ref{thm-main} and deduce some inequalities.

Before proving Theorem \ref{thm-main}, let us first give an overview of the method used here to prove the immanant positivity
of matrices, which was initiated by
Goulden and Jackson \cite{GJ92} and further developed by Greene \cite{Gre92},  Stembridge \cite{Ste92} and Wolfgang \cite{Wol97}. Goulden and Jackson \cite{GJ92} studied immanants of the Jacobi-Trudi matrices and conjectured that these immanants are nonnegative linear combinations of monomial symmetric functions. Based on a planar network interpretation for the Jacobi-Trudi matrices, they actually reduced their conjecture to a problem of characters of symmetric groups. To illustrate this, we need some notation. Let $J$ be a subinterval of $[n]=\{1,2,\ldots,n\}$. We denote by $\mathfrak{S}_J$ the subgroup of $\mathfrak{S}_n$ composed of all the permutations which fix all the elements in $[n]\setminus J$. Define $S_J = \sum_{\pi \in \mathfrak{S}_J} \pi$ as an element in the group algebra of $\mathfrak{S}_n$. We denote by $\Theta$ the set of all finite products of $S_J$'s. For convenience, we define the bracket of a square matrix $M=(m_{i,j})_{1\leq i,j\leq n}$ as
\[
    [M] = \sum_{\pi \in \mathfrak{S}_n} m_{1,\pi(1)}m_{2,\pi(2)}\cdots m_{n,\pi(n)} \cdot \pi.
\]
By naturally extending the characters of $\mathfrak{S}_n$ linearly to its group algebra,
one can easily find that $\mathrm{Imm}_{\lambda}M = \chi^{\lambda}([M])$ for any partition $\lambda$ of $n$.
Goulden and Jackson \cite{GJ92} proved that the coefficient of each monomial in the bracket of a Jacobi-Trudi matrix is a sum of elements in $\Theta$. Hence they reduced the problem of proving monomial nonnegativity to the following conjecture.

\begin{conj}[\cite{GJ92}]\label{conj-gj}
Let $\chi^{\lambda}$ be any irreducible character of $\mathfrak{S}_n$. Then for any $\theta \in \Theta$ the value $    \chi^{\lambda}(\theta)$ is nonnegative.
\end{conj}

Later, Greene \cite{Gre92} confirmed this conjecture by proving the following stronger result.
\begin{thm}[\cite{Gre92}]\label{thm-Greene}
Let $J$ be a subinterval of $[n]$, $\lambda$ be a partition of $n$ and $\rho^{\lambda}$ denote Young's seminormal representation of $\mathfrak{S}_n$ associated with $\lambda$. Then $\rho^{\lambda}(S_J)$ is a matrix with all entries being nonnegative.
\end{thm}

Based on Greene's result, Stembridge showed that the following result holds in the proof of \cite[Corollary 3.4]{Stembridge-1991}, which is actually equivalent to Conjecture \ref{conj-gj}.

\begin{prop}[{\cite{Stembridge-1991}}]\label{prop-stem}
Let $\chi^{\lambda}$ be any irreducible character of $\mathfrak{S}_n$. Then for any $\theta \in \Theta$,
\[
    \chi^{\lambda}(\theta)-\deg(\chi^{\lambda})\
    \chi^{(1^n)}(\theta) \ge 0,
\]
where $\deg(\chi^{\lambda})$ denotes the degree of $\chi^{\lambda}$.
\end{prop}

It should be noted that Haiman  proved a stronger result by using the Kazhdan-Lusztig theory, see \cite[Lemma 1.1]{Hai93}, and he further proved that the immanants of the Jacobi-Trudi matrices  are nonnegative linear combinations of Schur functions.

Stembridge \cite{Ste92} noticed that Goulden and Jackson's reduction procedure is carried out in the language of digraphs, but it uses special properties of the digraphs for Jacobi-Trudi matrices. Wolfgang \cite{Wol97} further studied general digraphs and obtained a reduction similar to that of Goulden and Jackson, which we will recall below.

Given a digraph $D$, we say that two paths intersect if they share a common vertex. Given two sequences of vertices $\mathbf{U} = (u_1,\ldots,u_n), \mathbf{V} = (v_1,\ldots,v_n)$ in $D$, we say that  $\mathbf{U}$ and $\mathbf{V}$ are  \textit{$D$-compatible} if any directed path from $u_{i}$ to $v_{j}$ intersects any directed path from $u_{k}$ to $v_{l}$ whenever $i < k$ and $j > l$.
Wolfgang \cite{Wol97} proved the following result, whose validity was surmised by Stembridge \cite[Section 6]{Ste92}.

\begin{thm}\cite[Section 2.6]{Wol97}\label{thm-d-compatible}
Let $D = (V(D),A(D))$ be a locally finite acyclic digraph, $\{z_a\mid a \in A(D)\}$ be a set of independent indeterminates and $(u_1,\ldots,u_n)$, $(v_1,\ldots,v_n)$ be two $D$-compatible sequences of vertices. Define the weight function $\mathrm{wt}_D$ of $D$ by $\mathrm{wt}_D(a) = z_a$ for each $a \in A(D)$. Then we have
\begin{equation*}
    \left[(GF_{D}(u_i,v_j))_{1\leq i,j\leq n}\right] = \sum_{\theta \in \Theta} f_{\theta} \cdot \theta,
\end{equation*}
where $f_{\theta}$ is a polynomial in $\{z_a\mid a \in A(D)\}$ with nonnegative coefficients,  $GF_{D}(u_i,v_j)$ denotes the sum of the weights of all directed paths from $u_i$ to $v_j$, and the weight of a directed path is the product of the weights of all its arcs.
\end{thm}

Applying Theorem \ref{thm-Greene} to Theorem \ref{thm-d-compatible}, one can easily obtain the following corollary.

\begin{cor}[{\cite[Corollary 2.6.5]{Wol97}}]\label{cor-wol97}
Let $D$, $z_a$, $(u_1,\ldots,u_n)$ and $(v_1,\ldots,v_n)$ be
as given in Theorem \ref{thm-d-compatible}.
Then for any irreducible character $\chi^{\lambda}$ of $\mathfrak{S}_n$, the immanant $\mathrm{Imm}_{\lambda}\left(GF_{D}(u_i,v_j)\right)_{1\leq i,j\leq n}$ is a polynomial in $\{z_a\mid a \in A(D)\}$ with nonnegative coefficients.
\end{cor}

We are now in the position to prove
Theorem \ref{thm-main}. Note that by suitably choosing sources and sinks, Wolfgang's theoretical framework actually applies to
any locally finite acyclic digraphs, and hence also applies to the planar networks constructed in Theorem \ref{thm-network-CS} and Corollary \ref{thm-network-CS-Hankel}.

\noindent \textit{Proof of Theorem \ref{thm-main}}.
We will only prove the immanant positivity of square submatrices of $C^{\gamma,\sigma,\tau}$. The proof for $H^{\gamma,\sigma,\tau}$ is similar and will be omitted.

Fix $n\geq 0$, let $I = (i_1,\ldots,i_m)$ and $J = (j_1,\ldots,j_m)$ be two sequences of indices such that $0 \le i_1<\cdots<i_m \le n$, $0 \le j_1<\cdots<j_m \le n$, and let $C_{I,J}$ be the submatrix of $C_n$ whose rows and columns are indexed by $I$ and $J$, respectively. Then
\[
\mathcal{C}_{I,J} = (D^{C_{I,J}},\mathrm{wt}_{D^{C_{I,J}}},(P_{n-i_1}^{(0)},\ldots,P_{n-i_m}^{(0)}), (P_{n-j_1}^{(n)},\ldots,P_{n-j_m}^{(n)}))
\]
is a planar network for $C_{I,J}$, where $D^{C_{I,J}}$ is
the subgraph of $D^{C_n}$ induced by the union of arcs of all directed paths from $P_{n-i_k}^{(0)}$ to $P_{n-j_l}^{(n)}$ for $1\leq k,l\leq m$, and $\mathrm{wt}_{D^{C_{I,J}}}$ is the restriction of $\mathrm{wt}_{D^{C_n}}$ to $D^{C_{I,J}}$. From our construction of $D^{C_n}$ in Section \ref{sect-ci}, it is easy to see that $(P_{n-i_1}^{(0)}, \ldots, P_{n-i_m}^{(0)})$ and $(P_{n-i_1}^{(n)}, \ldots, P_{n-j_m}^{(n)})$ are $D$-compatible in $D^{C_{I,J}}$.

Let $\{z_a\mid a \in A(D^{C_{I,J}})\}$ be a set of independent indeterminates. By Corollary \ref{cor-wol97}, for any character $\chi^{\lambda}$ of $\mathfrak{S}_m$, the immanant $\mathrm{Imm}_{\lambda} \ C_{I,J}$ is a polynomial in $z_a$'s with nonnegative coefficients. Then it suffices to show that $\mathrm{wt}_{D^{C_{I,J}}}$
assigns a $q$-nonnegative weight to each arc of $D^{C_{I,J}}$. For the first case of
Theorem \ref{thm-main}, we only need to take all $w_i$ in Theorem \ref{thm-network-CS} to be the weight function given by \eqref{eq-wt1}. The other four cases can be proved in the same manner. This completes the proof. \qed

We proceed to strengthen Theorem \ref{thm-main}. To this end, let us note the following stronger result than Corollary \ref{cor-wol97}, which can be obtained directly by applying Proposition \ref{prop-stem} to Theorem \ref{thm-d-compatible}.

\begin{cor} \label{cor-wol97-general}
Let $D$, $\{z_a\mid a \in A(D)\}$ and $(u_1,\ldots,u_n)$, $(v_1,\ldots,v_n)$ be
as given in Theorem \ref{thm-d-compatible}.
Then for any irreducible character $\chi^{\lambda}$ of $\mathfrak{S}_n$,
the difference
$$\mathrm{Imm}_{\lambda}\left(GF_{D}(u_i,v_j)\right)_{1\leq i,j\leq n}
-\deg(\chi^{\lambda})\cdot \det\left(GF_{D}(u_i,v_j)\right)_{1\leq i,j\leq n}$$ is a polynomial in $\{z_a\mid a \in A(D)\}$ with nonnegative coefficients.
\end{cor}

The next result follows from  Corollary \ref{cor-wol97-general} in the same way that Theorem \ref{thm-main} follows from Corollary \ref{cor-wol97}, so the proof is omitted.

\begin{cor}\label{cor-imm-ineq}
Let the three sequences $\gamma=(r_k)_{k\ge 0}$, $\sigma=(s_k)_{k \ge 0}$ and $\tau=(t_k)_{k\ge 1}$ be composed of $q$-nonnegative polynomials and satisfy one of the five conditions in Theorem \ref{thm-tp-csmatrix}.
Then for any $n\times n$ submatrix $M$ of $C^{\gamma,\sigma,\tau}$ or $H^{\gamma,\sigma,\tau}$, and for any irreducible character $\chi^\lambda$ of $\mathfrak{S}_n$, we have
$$\mathrm{Imm}_\lambda M - \mathrm{deg}(\chi^\lambda)\cdot \det M \ge_q 0.$$
\end{cor}

As mentioned in the introduction, one of the most important applications of the $q$-total positivity of $H^{\gamma,\sigma,\tau}$ is to prove the $q$-log-convexity of Catalan-like numbers. Since Corollary \ref{cor-imm-ineq}
establishes a stronger property of
$H^{\gamma,\sigma,\tau}$ than its $q$-total positivity, it is desirable to give more properties of the corresponding Catalan-like numbers. As an example of these applications, we have the following result.

\begin{cor}\label{cor-imm-3}
Let $\gamma$, $\sigma$ and $\tau$ be three sequences of $q$-nonnegative polynomials and satisfy one of the five conditions in Theorem \ref{thm-tp-csmatrix}, and denote $c_{k,0}$ by $a_k$ for any $k\geq 0$.
Then for any $0 \le i_1 < i_2 < i_3$ and $0 \le j_1 < j_2 < j_3$ we have
\begin{align}
2&(a_{i_1+j_2}a_{i_2+j_1}a_{i_3+j_3}
+a_{i_1+j_3}a_{i_2+j_2}a_{i_3+j_1}+a_{i_1+j_1}a_{i_2+j_3}a_{i_3+j_2}
)\nonumber\\
&-3(a_{i_1+j_2}a_{i_2+j_3}a_{i_3+j_1} +a_{i_1+j_3}a_{i_2+j_1}a_{i_3+j_2})\geq_q 0.\label{eq-application-1}
\end{align}
In particular, for any $0 \le i < j < k$ we have
\begin{align} \label{eq-application-2}
a_{2i}a_{j+k}^2
+a_{2j}a_{i+k}^2+ a_{2k}a_{i+j}^2 - 3a_{i+j} a_{j+k} a_{k+i}\geq_q 0.
\end{align}
\end{cor}

\pf It suffices to prove \eqref{eq-application-1}, and \eqref{eq-application-2} follows from \eqref{eq-application-1} by putting $i_1 = j_1 = i$, $i_2 = j_2 = j$, and $i_3 = j_3 = k$. Taking $n=3$ and $\lambda = (2,1)$ in Corollary \ref{cor-imm-ineq} and considering the submatrix of $H^{\gamma,\sigma,\tau}$ with row indices ${i_1,i_2,i_3}$ and column indices ${j_1,j_2,j_3}$ immediately lead to the desired result, by virtue of
\begin{align*}
&\deg(\chi^{(2,1)}) = \chi^{(2,1)}((1)(2)(3)) = 2, \\
&\chi^{(2,1)}((123)) = \chi^{(2,1)}((132)) = -1, \\
&\chi^{(2,1)}((23)(1)) = \chi^{(2,1)}((13)(2)) = \chi^{(2,1)}((12)(3)) = 0.
\end{align*}
This completes the proof. \qed

\noindent \textbf{Remark.} If $\gamma$, $\sigma$, and $\tau$ consist of real numbers, the inequality \eqref{eq-application-2} can be deduced
from the Hankel total positivity of $(a_n)_{n\geq 0}$
and the well-known Arithmetic-Geometric Mean Inequality. Precisely,
\begin{align*}
a_{2i}a_{j+k}^2
+a_{2j}a_{i+k}^2+ a_{2k}a_{i+j}^2 \ge & 3 \sqrt[3]{a_{i+j}^2a_{j+k}^2
a_{i+k}^2a_{2i}a_{2j} a_{2k}} \\
\ge & 3
\sqrt[3]{a_{i+j}^2a_{j+k}^2
a_{i+k}^2a_{i+j}a_{j+i} a_{2k}} \\
\ge & 3
\sqrt[3]{a_{i+j}^2a_{j+k}^2
a_{i+k}^2a_{i+j}a_{j+k} a_{k+i}} \\
= & 3a_{i+j} a_{j+k} a_{k+i}.
\end{align*}
However, we could not derive \eqref{eq-application-2} from the $q$-Hankel total positivity of $(a_n)_{n\geq 0}$
because the Arithmetic-Geometric Mean Inequality does not hold in general for polynomials.

Finally, we make use of Corollary \ref{cor-imm-3} to obtain some inequalities of some combinatorial sequences.

\begin{exam}
The Eulerian polynomials $E_n(q)$ \cite{StaEC1} are defined by
    \begin{equation*}
    \sum_{k\ge 0}(k+1)^nq^k=\frac{E_n(q)}{(1-q)^{n+1}},
    \end{equation*}
     which can be seen as Catalan-like numbers generated by the  Catalan-Stieltjes matrix $C^{\gamma,\sigma,\tau}$
     with $r_k=k+1$, $s_k=k(q+1)+1$ and $t_{k+1}=(k+1)q$ for all $k\ge 0$, see \cite{Pan-Zeng-2016}. It is clear that $\gamma$, $\sigma$, and $\tau$ satisfy the first condition of Theorem \ref{thm-tp-csmatrix}. By Corollary \ref{cor-imm-3}, for any $0\leq i<j<k$ we have
\begin{align*}
E_{2i}(q)E_{j+k}^2(q)
+E_{2j}(q)E_{i+k}^2(q)+ E_{2k}(q)E_{i+j}^2(q) - 3E_{i+j}(q) E_{j+k}(q) E_{k+i}(q)\geq_q 0.
\end{align*}

\end{exam}

\begin{exam}
The Schr\"{o}der polynomials $R_n(q)$ \cite{Bonin-1993} are defined by
    \begin{equation*}
    R_n(q)=\sum_{k=0}^{n}\binom{n+k}{n-k}\frac{1}{k+1}\binom{2k}{k}q^{k}.
    \end{equation*}
    They are Catalan-like numbers generated by the  Catalan-Stieltjes matrix $C^{\gamma,\sigma,\tau}$
     with $r_k=1$ for $k\ge 0$, $s_0=q+1$, $s_k=2q+1$ and $t_{k}=q(q+1)$ for $k\ge 1$, see  \cite{Wang-Zhu-2016,Zhu-2013}.
     Obviously, $\gamma$, $\sigma$, and $\tau$ satisfy the fifth condition of Theorem \ref{thm-tp-csmatrix}.
     It follows from Corollary \ref{cor-imm-3} that for any $0\leq i<j<k$,
\begin{align*}
R_{2i}(q)R_{j+k}^2(q)
+R_{2j}(q)R_{i+k}^2(q)+ R_{2k}(q)R_{i+j}^2(q) - 3R_{i+j}(q) R_{j+k}(q) R_{k+i}(q)\geq_q 0.
\end{align*}
\end{exam}

\begin{exam}
The Narayana polynomials of type A \cite{Petersen2015Eulerian} are defined by
    \begin{equation*}
    N_n(q)=\sum_{k=1}^{n}\frac{1}{n}\binom{n}{k-1}\binom{n}{k}q^k
    \end{equation*}
    for $n \ge 1$ and $N_0(q) := 1$.
    They are Catalan-like numbers generated by the  Catalan-Stieltjes matrix $C^{\gamma,\sigma,\tau}$
     with $r_k=1$ for $k\ge 0$, $s_0=q$, $s_k=q+1$ and $t_{k}=q$ for $k\ge 1$, see \cite{Wang-Zhu-2016,Zhu-2013}.
     It is straightforward to verify that $\gamma$, $\sigma$, and $\tau$ satisfy the second, fourth, and fifth conditions of Theorem \ref{thm-tp-csmatrix}.
     Then by Corollary \ref{cor-imm-3}, for $0 \le i < j < k$ we have
\begin{align*}
N_{2i}(q)N_{j+k}^2(q)
+N_{2j}(q)N_{i+k}^2(q)+ N_{2k}(q)N_{i+j}^2(q) - 3N_{i+j}(q) N_{j+k}(q) N_{k+i}(q)\geq_q 0.
\end{align*}
\end{exam}

\noindent \textbf{Acknowledgments.} This work is supported in part by the Fundamental Research Funds for the Central Universities and the National Science Foundation of China (Nos. 11522110, 11971249).

\end{document}